\documentclass[twoside,reqno]{lt6-proc}
\usepackage{epsfig,cite}
\usepackage{amssymb,amsmath}
\usepackage{times}
\begin{document}
\sloppy \raggedbottom
\setcounter{page}{1}

\def\oplusinf{\mathop{\oplus}}
\def\otimesinf{\mathop{\otimes}}

\def\id{{\mbox{Id}}}
\def\im{{\mbox{Im}}}
\def\ker{{\mbox{Ker}}}
\def\Hom {{\mbox{Hom}}}
\def\hom {{\mbox{hom}}}
\def\End{{\mbox{End}}}
\def\Ext{{\mbox{Ext}}}
\def\cali{{\cal I}}
\def\calp{{\cal P}}
\def\cala{{\cal A}}
\def\calb{{\cal B}}
\def\cald{{\cal D}}
\def\calh{{\cal H}}
\def\cals{{\cal S}}
\def\calk{{\cal K}}
\def\fraka{\mathfrak A}
\def\frakc{\mathfrak C}
\def\alg{{\mathbf{Alg}}}
\def\vect{{\mathbf{Vect}}}
\def\Homg {{\mathbf{Hom}}}
\def\homg {{\mathbf{hom}}}
\def\endg{{\mathbf{end}}}
\def \la {\left\langle}
\def \ra {\right\rangle}

\newcommand{\vir}{\raisebox{0.75mm}{,}}
\newtheorem{theorem}{Theorem}
\newtheorem{lemma}{Lemma}
\newtheorem{proposition}{Proposition}
\newtheorem{corollary}{Corollary}
\newtheorem{definition}{Definition}
\newcommand{\beq}{\begin{equation}}
\newcommand{\eeq}{\end{equation}}
\newcommand{\beqa}{\begin{eqnarray}}
\newcommand{\eeqa}{\end{eqnarray}}
\newcommand{\beqann}{\begin{eqnarray*}}
\newcommand{\eeqann}{\end{eqnarray*}}
\newcommand{\nn}{\nonumber \\}
\newpage
\setcounter{figure}{0}
\setcounter{equation}{0}
\setcounter{footnote}{0}
\setcounter{table}{0}
\setcounter{section}{0}



\title{Automorphisms of Regular  Algebras}

\runningheads{Todor Popov}{Automorphisms of Regular Algebras}

\begin{start}

\author{Todor Popov}{1}

\address{Institute for Nuclear Research and Nuclear Energy,\\
Bulgarian Academy of Sciences, \\ 
bld. Tsarigradsko chaus\'ee 72,
Sofia 1784}{1}

\begin{Abstract}
Manin associated to a quadratic algebra (quantum space) the quantum matrix 
  group of its automorphisms. This Talk aims to demonstrate that Manin's construction 
can be extended for quantum spaces which are non-quadratic homogeneous algebras. 
Here given  a  regular Artin-Schelter  algebra  of 
dimension 3 we construct  the quantum group of its  symmetries, i.e.,
the Hopf algebra of its automorphisms.
For quadratic Artin-Schelter algebras these quantum groups are contained in the
the classification of the  $GL(3)$ quantum matrix groups due to Ewen and Ogievetsky.
For   cubic  Artin-Schelter  algebras we obtain new quantum groups
which are  automorphisms of  cubic quantum  spaces. 
\end{Abstract}
\end{start}


All vector spaces and algebras are  over a ground field   $\mathbb K$
of characteristics $0$. We adopt the Einstein  convention of summing
on repeated an upper and a lower indices  except when these  are in
brackets, e.g. there is no summation in $Q^{(i)}_{(i)}$.\\

\section{$N$-Homogeneous Algebras}

A  $N$-homogeneous algebra is an algebra of the form \cite{RB1}, \cite{BDVW}
$$\cala=A(E,R)=T(E)/(R)$$
where $E$ is finite dimensional vector space, $T(E)$ is the tensor 
algebra of $E$ and $(R)$ is the two-sided ideal generated by a vector 
subspace $R\subset E^{\otimes N}$.  Since the space $R$ is homogeneous
by ascribing the degree 1 to the generators in $E$ 
one obtains that the algebra $\cala$ is graded, $\cala=\oplus_{n\geq 0} \cala_{n}$,
generated in degree 1,  $\cala_{0}=\mathbb K$ and such that 
the degrees $\cala_{n}$ are finite-dimensional vector spaces.

The dual  $\cala^{!}$ of $\cala=A(E,R)$ is defined to be the
$N$-homogeneous algebra $\cala^{!}=A(E^{\ast},R^{\perp})$
where $E^\ast$ is the dual vector space of $E$ and 
$R^{\perp}\subset E^{\ast \otimes N}= (E^{\otimes N})^{\ast}$
is the annhilator of $R$, $R^{\perp}(R)=0$ . One has $(\cala^{!})^{!}=\cala$. \\
  
\noindent
 Given two $N$-homogeneous algebras $\cala=A(E,R)$ and
 $\cala'=A(E',R')$   one defines the $N$-homogeneous algebra
 $$\cala \bullet \cala'=A(E\otimes E', \pi_{N}(R\otimes R'))$$
where $\pi_{N}$ is the permutation $$\pi_{N}(e_1\otimes \ldots \otimes e_N
   \otimes e_1^{'} \otimes \ldots \otimes e_N^{'} )=
    e_1\otimes e_1^{'} \ldots e_N\otimes e_N^{'}.$$

    Following Manin's monograph\cite{YM1} on quadratic 
 algebras  R. Berger, M. Dubois-Violette and M. Wambst\cite{BDVW} 
 introduced   the corresponding semigroup $end(\cala)=\cala^{!}\bullet \cala$
 of the endomorphisms of the  $N$-homogeneous algebra $\cala$.
The semigroup  $end(\cala)$ is canonically endowed with the structure of  bialgebra 
with a coproduct $\Delta$ and counit $\varepsilon$ given by
   \beq
   \Delta (u_{j}^{i}) = u_{k}^{i}\otimes u_{j}^{k}
   \qquad 
   \varepsilon(u_{i}^{j}) = \delta_{i}^{j}.
   \label{bialg}
   \eeq 
The algebra $\cala$ is a left comodule of $end(\cala)=\cala^{!}\bullet \cala$
 for the coaction\cite{BDVW} 
 \beq
 \delta(x^{i}) = u_j^{i} \otimes x^{j} 
 \qquad u_j^{i} \in end(\cala) =A(E^{\ast}\otimes E, r)
 \qquad r=\pi_{N}(R^{\perp}\otimes R).
   \eeq  

 The semigroup  $end(\cala)$   alone has not enough relations in order to allow 
for an antipode. In order to obtain   the quantum matrix group(the group 
of the automorphisms of  $\cala$) we shall proceed by ``adding the missing 
relations''\cite{YM1}, considering also the semigroup $end(\cala^{!})$ of the endomorphisms of the dual  $\cala^!$.

\noindent
The algebra  $\cala^{!}$ becomes a left comodule of the semigroup
 $end(\cala^{!})=\cala\bullet \cala^{!}$ in the  following way. 
 Let us identify the upper and the lower indices in $\cala^{!}$
with the help of the bilinear form  $g^{ij}=\delta^{ij}$, 
$
\xi^i= g^{ij} \xi_{j} =\xi_{i}
$.
A left coaction on $\cala^{!}$ is given by
 \beq
   \delta(\xi^{i}) = \check{u}_j^{i} \otimes \xi^{j}
   \quad \check{u}_j^{i} \in end(\cala^{!})_{g} =A(E^{\ast}\otimes E,\check{r})
   \qquad \check{r}=\pi_{N}(R\otimes R^{\perp})^{\ast}
   \label{coa}
   \eeq
   where $end(\cala^{!})_{g}$ is $end(\cala^{!})$ up to a move of the 
   indices (by $g$). 
   We  shall not distinguish $end(\cala^{!})_{g}$ and $end(\cala^{!})$,
   thus $\cala^{!}$ is a left   $end(\cala^{!})$-comodule.   
   
   Let us consider the bialgebra $e(\cala)$(again in the spirit of  \cite{YM1}) having 
   relations those of the semigroups $end(\cala)$  and $end(\cala^{!})$,
   in which the generators  ${u}_j^{i}$, $\check{u}_j^{i}$ 
   have been identified ${u}_j^{i} \equiv\check{u}_j^{i}$, i.e., we 
   consider the bialgebra   
  \beq
  e(\cala)= A(E^{\ast}\otimes E,r\oplus \check{r})/({u}_j^{i} -\check{u}_j^{i})
 \qquad \Delta (u_{j}^{i}) = u_{k}^{i}\otimes u_{j}^{k}
   \qquad 
   \varepsilon(u_{i}^{j}) = \delta_{i}^{j}
  \eeq
which is quotient of the bialgebra $end(\cala)$ (and $end(\cala^{!})$). 
The algebras $\cala$ and $\cala^{!}$ are left $e(\cala)$-comodules 
in a natural way.

\section{Regular Artin-Schelter Algebras}
The question that we address in the Talk is
 when an algebra $\cala$ is ``good'' quantum space, in the sence  that
 the space of its endomorphism  $e(\cala)$ is a Hopf algebra, i.e.,
there is a  quantum  group of the symmetries of the quantum space $\cala$?

Artin and Schelter have considered a class of regular algebras with 
very ``good'' homological properties \cite{AS}(we use here the equivalent 
definition of \cite{YM2}).

\begin{definition}
A graded algebra $A=\bigoplus_{n\geq 0} A_{n}$ with $A_{0}=\mathbb K$, 
generated by $A_{1}$, $\dim A_{1}< \infty $ is called regular
    if:
    
i) $A$ has polynomial growth, 
(i.e., $gk$-$\dim A = \gamma < \infty $),
    
ii) it is Gorenstein, i.e., there is a finite free 
resolution of the trivial right $A$-module ${\mathbb K}$, such that
its dualized complex (by $Hom_{A}(\bullet, A)$)  is a finite free resolution 
of the trivial left $A$-module ${\mathbb K}$. The length of this 
resolution is called dimension of $A$.

 %
 %
    \end{definition}

    Manin suggested \cite{YM2} 
    that the regular algebras are  good candidates for  quantum spaces. 
The regular algebras of dimension 2
are exhausted by the Manin plane $yx-qxy=0$  and the Jordanian plane $xy-yx-y^2=0$
(given in the Introducion of \cite{AS} as  simple examples). 

The classification of  the regular algebras $A$ of 
 dimension 3 which is done in \cite{AS} is much more involved and 
 requires some new technics;
it turns out that a regular algebra $A$ of dimension 3 is 
 either quadratic ($N=2$) or cubic ($N=3$) homogeneous algebra,
 $\cala=A(E,R)$, generated by two elements satisfying two cubic relations,
 or else by three elements with three quadratic relations.
 Further the classification is based on the fact that a regular algebra $\cala$
 has an intrinsic description in terms of an invariant element
 $\omega(\cala)$ of degree $N+1$ in the generators,
 i.e., in terms of a tensor with $N+1$ indices(Proposition (2.4) of \cite{AS}).
The invariant tensor $\omega(\cala)$ transforms as the coordinates of
 1-dimensional space  of the 
the maximal  non-vanishing degree $\cala^{!}_{m}$  of $\cala^{!}$
\beq
    \cala^{!}_{m} = E^{\ast}\otimes
   R^{\ast} \cap R^{\ast} \otimes E^{\ast}, \qquad m=N+1.
   \label{erre}
   \eeq
Thus the invariant $\omega$ of the algebra $\cala$ can be written 
in  either of the ways
   \beq
  \omega =\xi^{i}f_i^\ast = g_{i}^{\ast} \xi^{i}= 
Q_{i}^{j}f_{j}^{\ast}
  \xi^{i}
  \eeq
  where $f_{j}$ and $g_{i}=Q_{i}^{j} f_{j}$ are two bases of $R$.
By change of the basis the matrix $Q$ is amenable to 
  the Jordanian canonical form. Artin and Schelter 
  have proven that the case of $Q$ diagonal matrix  is generic,
  in the sense that the non-diagonal matrices $Q$ do not 
  give new regular algebras.  On the  components of  
   $\omega(\cala)=\omega_{I}\xi^{I} \in \cala^{!}_{N+1}$ 
 is defined the cyclic action 
  \beq
 \sigma(\omega_{Ai})=\omega_{i A} = Q_{i}^i\omega_{Ai} \qquad
 \sigma(\omega_{i_1 \dots i_{N} i_{N+1}}) = \omega_{ i_{N+1} i_1 
\dots
   i_{N} }= \omega_{\sigma^{-1}(I)},
   \label{cycl}
   \eeq 
   under which $\omega$ splits into orbits.
    If the multiindex $J=j_1 \dots j_n$ belongs to  an orbit of the cyclic action  of order $n$, i.e., if $\sigma^n(J)=J$ 
   then  $\omega_{J}= 0$ or else $ Q_{j_{1} }^{j_{1}} \ldots Q_{j_{n} }^{j_{n}} =1$. 
 
 An algebra $\cala$ in the classication of the regular algebras of 
 dimension 3 is characterised by 
 the data $(Q,\omega)$, the diagonal matrix $Q$ and the invariant $\omega$
 (see Table (3.9) and (3.11) in \cite{AS})
where  are listed  6 families of cubic algebras 
(types $A$,$E$,$H$,$S_{1}$,$S_{2}$,$S_{2}'$)
and 7 types of quadratic algebras(types 
$A$,$B$,$E$,$H$,$S_{1}$,$S_{1}'$,$S_{2}$).

The cubic class $S_1$ contains  the universal 
enveloping algebra of the Heisenberg algebra \cite{AS} with two
generators (then coinciding with the Yang-Mills and related algebras \cite{AC.MDV1},
\cite{AC.MDV2}) as well as the algebras related to the parastatistics \cite{D-VP}. 

\section{Koszul Complexes for $N$-Homogeneous Algebras}
To every $N$-homogeneous algebra $\cala$ one associates canonically 
the {Koszul  complexes } ${\cal L}(\cala)$ and ${\cal K}(\cala)$ 
   \cite{BDVW}, \cite{BM}(see also \cite{AC.MDV1}, \cite{AC.MDV2})
   as follows.

Let us introduce the element $c$, independent of the choice of  basis
$$c=  \xi_{i} \otimes x^i\in  \cala^! \otimes \cala.$$ 
and the linear mapping $d$ defined by the left multiplication of $c$,
$$ d:\cala^!_n \otimes \cala \rightarrow \cala^!_{n+1} \otimes \cala \qquad d: \alpha \otimes a \mapsto c(\alpha \otimes a) = 
\xi_{i} \alpha  \otimes x^ia.$$
In view of the definition of $R^{\perp}$ the canonical element satisfies $c^N=0$, which implies $d^N=0$. In other words the mapping $d$ is a $N$-differential mapping and the sequence of spaces $(\oplus_{n\geq 0}\cala^!_n \otimes \cala ,d)$ is a $N$-complex (for details we send the reader to \cite{BDVW}, \cite{BM}).

The {\it Koszul cochain complex} $({\cal L}(\cala),{\mathfrak d})$ 
is defined to be the complex \cite{BDVW}
$$({\cal L}(\cala),{\mathfrak d})=
(\bigoplus_{i\geq 0} {\cal L}^{i}(\cala),\bigoplus_{i\geq 0} {\mathfrak d}^i)$$
with 
 degrees given by
\beq
{\cal L}^{2i} (\cala) =  \cala^{!}_{Ni} \otimes \cala 
\qquad {\cal L}^{2i+1} (\cala) =  \cala^{! }_{Ni+1} \otimes \cala 
\label{coKosc}
\eeq
and  differential mapping 
${\mathfrak d}^j: {\cal L}^{j} (\cala)\rightarrow {\cal L}^{j+1} (\cala)  $, 
 $$ 
 {\mathfrak d}^j = \left\{\begin{array}{lcl} d &\mbox{when}& j=2i,\\
  d^{N-1} &\mbox{when}& j=2i+1. \end{array} \right. 
 $$
 The jump in the degrees (when $N>2$)
${\cal L}^i(\cala)=
\cala^{!}_{n(i)} \otimes \cala $, 
$n(2i)=Ni$ and $n(2i+1)=Ni+1$, 
is due to the fact that the complex ${\cal L}(\cala)$ is a contraction of an  underlying $N$-complex \cite{BDVW}.

The {\it Koszul chain complex} $({\cal K}(\cala),{\mathfrak d}')$ \cite{BDVW} can be defined as the dualized  complex of the cochain complex ${\cal L}(\cala)$, 
${\cal K}(\cala)=Hom_{\cala}({\cal L}(\cala), \cala)$ with degrees
\beq
{\cal K}(\cala)=\bigoplus_{i\geq 0} {\cal K}_{i}(\cala)=\bigoplus_{i\geq 0}\cala \otimes \cala_{n(i)}^{!\ast}
\label{Kosc}\eeq
and differential ${\mathfrak d}'$ defined by $d'$ on ${\cal K}_{2i}$
and $d'^{N-1}$ on  ${\cal K}_{2i+1}$ where $d' $
is the mapping 
 $$d':\cala \otimes \cala_{k}^{!\ast} \rightarrow \cala \otimes \cala_{k-1}^{!\ast}
 \qquad d': x_0 \otimes (x_1\otimes x_2 \otimes \dots x_n) \mapsto x_0  x_1\otimes(x_2 \otimes \dots x_n)$$
 for which $d'^N=0$ holds.
The dualization (by $Hom_{\cala}(\bullet, \cala)$) of the Koszul chain complex  ${\cal K}(\cala)$  gives back the cochain complex, ${\cal L}(\cala)=Hom_{\cala}({\cal K}(\cala), \cala)$.

It was shown in \cite{BDVW} (using contraction of $N$-complexes)
that the notion of Koszul algebra has a meaningful generalization
 for 
$N$-homogeneous algebras;
a $N$-homogeneous algebra $\cala$ is said to be {\it Koszul algebra} when its Koszul chain complex ${\cal K}(\cala)$ 
is acyclic in positive degrees \cite{BDVW}. 
 The  Koszul chain complex ${\cal K}(\cala)$ of a Koszul algebra
$\cala$ provides a free resolution of the trivial left module $\mathbb K$,
property which we shall use right now.

\section{ Homological Quasi-Determinant}
We are ready to introduce the analog of the determinant for a regular Koszul $N$-homogeneous  
algebra $\cala$. 
 Let $\cala$ be $N$-homogeneous Koszul and regular of 
dimension $D$. Then the Koszul chain complex ${\cal K}(\cala)$ (\ref{Kosc}) of $\cala$ provides a free resolution of the left $\cala$-module
$\mathbb K$ with length $D$ and the Gorenstein property implies that the  complex  ${\cal L}(\cala)$ (\ref{coKosc}) provides a resolution of the trivial right $\cala$-module $\mathbb K$,
i.e., $$H({\cal L}(\cala))= H^{D}({\cal L}(\cala))= \cala^{!}_{n(D)}=\cala^{!}_{max}\simeq  \mathbb K,$$
thus the cohomology group $\cala^{!}_{max}$ of ${\cal L}(\cala)$ is a $1$-dimensional 
$e(\cala)$-comodule.
\begin{definition}
The bialgebra $e(\cala)$ of a regular Koszul $N$-homogeneous   algebra $\cala$
coacts on $\cala^{!}_{max}$ 
by the element  $\cald=\cald(\cala)$ 
  \beq
    \delta:\cala^{!}_{max}\rightarrow e(\cala)\otimes \cala^{!}_{max}, 
\qquad
    \delta( \omega(\cala)) = \cald(\cala) \otimes \omega(\cala)
    \label{qdet}
    \eeq  
  which we are referring to as the {\it homological quasi-determinant}. 
  \end{definition}
   For the quadratic algebras $N=2$,  the homological quasi-determinant 
   $\cald$  coincides  with the  Manin's { \it homological determinant} \cite{YM1}. 
  
\begin{lemma}
      Let 
      us 
      denote the one-dimensional comodule $\cala^{!}_{max}$ as $\omega(\cala)=\omega_{A}\xi^{A}:= \omega_{a_1 \ldots
        a_{m}}\xi^{a_1} \ldots \xi^{a_{m}} $.
        If the coefficient $\kappa=\omega^{A}\omega_{A}$ 
	 is invertible
       then the homological quasi-determinant $\cald$ of the bialgebra 
       $(e(\cala), \Delta, \varepsilon)$ reads
           \beq
        \cald = \kappa^{-1} \omega_{A} u^{A}_{B}  \omega^{B}=
         \kappa^{-1} \omega_{a_1 \ldots a_m}
         u^{a_1}_{b_1} \ldots u^{a_m}_{b_m}  \omega^{b_1 \ldots b_m}
         \label{cqdet}
    \eeq
    The element $\cald$ of $e(\cala)$
    is group-like, $\Delta (\cald)=\cald  \otimes \cald$ and $\varepsilon(\cald)=1$.
    \label{quasidet}
\end{lemma}
{\it Proof:} In components the coaction  
of the bialgebra $e(\cala)$  on the $1$-dimensional comodule
$\omega \in \cala_{m}^{!}$ (\ref{qdet}) reads
\beq
\omega_{A} u^{A}_{B} =\cald \omega_{B}.
\eeq
Multiplying of both sides by $\omega^{B}$ and summing on the 
multi-index $B$ yields
\beq
\omega_{A} u^{A}_{B}\omega^{B}=\cald \omega_{B}
\omega^{B} = \kappa \cald
\eeq
which for invertible $\kappa$ implies the quasi-determinant formula (\ref{cqdet}). 
The coproduct and the counit on a monomial $u^A_{B}$ is given by
$\Delta(u^A_{B})=u^A_{C}\otimes u^C_{B}$   and $\varepsilon(u^{A}_{B})=\delta^{A}_{B}$
hence
\beq
\begin{array}{rcccccl}
\Delta (\cald)  &=&
\kappa^{-1} \omega_{A}  u^{A}_{C}\otimes u^{C}_{B} \omega^{B}&=&
\kappa^{-1} \cald \omega_{C} \otimes u^{C}_{B} \omega^{B} &=& \cald
\otimes \cald, \\
\varepsilon(\cald) &=& \kappa^{-1} \omega_{A} \varepsilon(u^{A}_{B})
\omega^{B} &=&
\kappa^{-1} \omega_{A} \delta^{A}_{B}  \omega^{B}
&=&1. \qquad \Box
\end{array}
\label{eps}
\nonumber
\eeq

For the quasi-determinant $\cald$ there exist an expansion 
on rows and on columnes which gives rise to the Cramer adjoint elements
(analogues of  subdeterminants).
 \begin{definition}
  The left $\cals_{L}(u_{j}^{k})$ and right $\cals_{R}(u_{j}^{k})$
  Cramer adjoint elements 
  are such elements of the bialgebra $(e(\cala),\Delta, \varepsilon)$ 
  that the following holds
  \beq
   \cals_{L}(u_{k}^{i}) u_{j}^{k}=\cald  \delta_j^i = \cald
   \varepsilon( u_{j}^{i}) 
    =u_{k}^{i} \cals_{R}(u_{j}^{k}) .
   \label{cram}
    \eeq
  \end{definition}

\section{Automorphisms of Regular Algebras of Dimension 3}   
 The regular algebras of dimension 2, i.e., the Manin and the Jordanian plane are Koszul algebras.
 The quadratic and cubic 
  regular algebras of dimension  3 \cite{AS} are also Koszul algebras \cite{RB1},\cite{BM} which allows to define their homological quasi-determinants $\cald$  coinciding with the usual quantum determinants \cite{YM1} for $N=2$, but giving something new for cubic algebras $N=3$.
 
 From now on $\cala$ will stay for a regular  (quadratic or cubic) algebra of  dimension 3 of the Artin-Schelter clasification,  
specified by its  data $(Q,\omega)$ \cite{AS}. We put the accent
on the cubic Artin-Schelter algebras, but the quadratic ones fit into the same description, which allows for their simultaneous treatment.
  
In our construction of the antipode on the bialgebra $e(\cala)$
 the Cramer 
  adjoint elements will be instrumental. 
\begin{definition}
    A regular algebra $\cala$ of dimension 3 will be  referred to as  generic regular
    algebra  when the coefficients  $\kappa^{i}$ 
    ($\tilde{\kappa}_{i}$)
    are invertible
 \beq
    {\kappa}^{i}= \omega^{A (i)} \omega_{A (i)}  \neq 0
    \qquad (\tilde{\kappa}_{i} = \omega^{(i)A}\omega_{(i)A} \neq 0) 
    \quad \mbox{no summation on $(i)$!}
    \label{coef}
    \nn
    \eeq
The bialgebra $e(\cala)$ for a generic $\cala$ will  be referred 
to as generic bialgebra.
    \end{definition}
The equation ${\kappa}^{i}=0$ for some $i$ (equivalent to 
$\tilde{\kappa}_{i}=0$  in view of the cyclicity (\ref{cycl}))  
singles out the non-generic algebras. 

    One easily checks that
the  regular cubic  algebras $\cala$ of type $E$,$H$ and
some points of the types $A$, $S_{1}$ and  $S_{2}$  are not generic.

    \begin{proposition} Let $\cala$ be a generic regular algebra
        of dimension 3.
  The left(right) Cramer adjoint elements of  the bialgebra
  $e(\cala)$  are given by
   \beq
    \cals_{L}(u^{j}_{i}) = ({{\kappa}^{j}})^{-1} \omega_{ {A} i} 
u^{{A}}_{{B}}
\omega^{{B} j}   \qquad
(\cals_{R}(u^{j}_{i}) = ( \, \tilde{\kappa}_{i})^{-1}
\omega_{i {A}} u_{ {B}}^{{A}} \omega^{j {B}} \,).
  \label{cramer}
\eeq
 \end{proposition}
{\it Proof:} For a generic  bialgebra $e(\cala)$  one has
$\omega^{ {A} i} \omega_{{A} j}={\kappa^{i}}\delta_{j}^{i}$
with $\kappa^{i} \neq 0$ hence
    \beq
\cals_{L}(u_{k}^{i}) u_{j}^{k}=({{\kappa}^{i}})^{-1}
\omega_{ {A} k} u^{ {A}k}_{{B}j} \omega^{{B} i }=
\cald \omega_{ {B} j} \omega^{{B} i} ({{\kappa}^{i}})^{-1}=
\cald  \delta_j^i = \cald \varepsilon(u_{j}^{i}).
\nonumber
\eeq
The right Cramer adjoint element are handled similarly due to
$\omega^{ i{A} } \omega_{j {A} }=\tilde{\kappa}_{j}\delta_{j}^{i}$. 
$\Box$


\begin{proposition} Let $\cala$ be a generic regular algebra
        of dimension 3. The left $\cals_{L}(u_{j}^{i})$ and the right
    $\cals_{L}(u_{j}^{i})$ Cramer adjoint elements in the bialgebra 
    $e(\cala)$ 
    are proportional 
     \beq
\cals_{L}(u_{j}^{i}) = h_{(j)}^{(i)}\cals_{R}(u_{j}^{i}) \qquad 
\eeq
with coefficient a constant $h_{(j)}^{(i)}$ 
being an element of the multiplicatively antisymmetric matrix
\beq
h_{(j)}^{(i)}= \left(\frac{\tilde{\kappa}_{(i)}Q^{(i)}_{(i)}}
{\tilde{\kappa}_{(j)}Q^{(j)}_{(j)}}\right)^{-1}= \frac{1}{h_{(i)}^{(j)}},
\quad \mbox{}
\quad 
h_{(i)}^{(i)}=1, \qquad h_{(j)}^{(i)}h_{(k)}^{(j)}=h_{(k)}^{(i)}.
\label{h}
\eeq
\end{proposition}
{\it Proof:} Taking into account the 
cyclicity
$
\sigma(\kappa^{i})=  \tilde{\kappa}_{i} = (Q^{(i)}_{(i)})^2 
\kappa^{(i)}
$
one can bring the expression (\ref{cramer}) for the 
$\cals_{L}(u_{j}^{i})$  to the form of $\cals_{R}(u_{j}^{i})$ which 
gives the result. $\Box$

For the generic regular cubic  algebra $\cala$ which are our main 
concern here the formula (\ref{h}) yields 
$$
h_{(j)}^{(i)} =h^{j-i}$$
with $h=1$ for $\cala$ of cubic types $A$ and $S_{1}$, $h=-1$ for $S_{2}$
and  $h=-2/3$ for $S_{2}'$.

\begin{corollary}  The quasi-determinant $\cald=\cald(\cala)$ of the
 bialgebra $(e(\cala),\Delta, \varepsilon)$ of 
 a generic regular algebra $\cala$ is a quasi-central element
     \beq
\cald u_{j}^{i}  = h^{(i)}_{(j)} u_{j}^{i}\cald.
\eeq
\end{corollary}
{\it Proof:} By expansion of 
$\cald$ on the RHS  and regrouping the terms we get
\beqa
h^{(i)}_{(j)} u_{j}^{i}\cald & = &h^{(i)}_{(k)} u_{k}^{i}\cald 
\delta^k_{j}=
h^{(i)}_{(k)} u_{k}^{i} \cals_{L}(u_{s}^{k})u_{j}^{s}= 
h^{(i)}_{(k)} h^{(k)}_{(s)} u_{k}^{i} 
\cals_{R}(u_{s}^{k})u_{j}^{s} \nn &=& h^{(i)}_{(s)} u_{k}^{i} 
\cals_{R}(u_{s}^{k})u_{j}^{s}= h^{(i)}_{(s)}  \cald 
\delta^i_{s}u_{j}^{s} = h^{(i)}_{(i)} \cald u_{j}^{i} =  \cald u_{j}^{i}.
\nonumber
\eeqa
Note that the cancelation of the index $k$ in $h^{(i)}_{(k)} h^{(k)}_{(s)}=h^{(i)}_{(s)}$ 
is crucial for resummation of the terms. $\Box$
\begin{theorem}
       Let $\cala$ be a generic regular 
       algebra of 
       dimension 
       $3$\cite{AS}.
       Let us denote by $(\calh(\cala), \Delta, \varepsilon )$
       the 
       bialgebra  $(e(\cala), \Delta, \varepsilon )$ 
       extended
      by the inverse element $\cald^{-1}$, 
      $\cald \cald^{-1}=\cald^{-1}\cald={1\!\!1}_{e(\cala)}$
       and consider 
       the linear antihomomorphism 
       \beq
  S: \calh(\cala) \rightarrow 
       \calh(\cala)^{op} \qquad S: u_{i}^{j} \mapsto  {S}(u_{i}^{j})= \cald^{-1}\cals_{L} (u_{i}^{j})=
      \cals_{R}(u_{i}^{j})\cald^{-1}. 
      \label{anti}
      \nonumber
      \eeq
       Then  $(\calh(\cala), \Delta, \varepsilon, S)$ is a Hopf 
       algebra with an antipode given by $S$.   
         \nonumber
     \end{theorem}

     {\it Proof:} The bialgeba structure of $e(\cala)$ is compatible with the antipode,
     if they satisfy the antipode axiom 
\beq
m\circ(\id\otimes S )\circ\Delta =
m\circ( S\otimes \id )\circ \Delta = \eta\circ \varepsilon
\label{ax}
\eeq
where $m$ 
is the product  and $\eta$ is the unity mapping of the algebra $e(\cala)$,
$\eta: 
1 \mapsto {1\!\!1}_{e(\cala)} $.
The existence of the left $\cals_{L}$ and right $\cals_{R}$ Cramer adjoint 
     elements for  $e(\cala)$ of a generic algebra $\cala$(Proposition 1) 
    implies that the antipode  $S$ constructed by $\cals_{L}$ 
    or $\cals_{R}$ 
    satisfies the axiom (\ref{ax}) which makes the bialgebra $e(\cala)$ a Hopf algebra.$\Box$\\

\noindent        
 The antipode $S(\cald)$ of the quasi-determinant $\cald$ 
 is evaluated by 
the axiom (\ref{ax}) 
\beq
S(\cald)\cald=\cald S(\cald)= \eta\circ
\varepsilon(\cald)=  1\!\!1_{e(\cala)} 
\nonumber
\eeq
where we used $\Delta(\cald)=\Delta(\cald)\otimes \Delta(\cald)$
and $\varepsilon(\cald)=1$ (Lemma \ref{quasidet}). 
Thus ${S}(\cald)=\cald^{-1}$. \\

\section{Conclusions}
We have constructed the Hopf algebra $\calh(\cala)$ of the 
automorphisms of a generic regular algebra $\cala$ of dimension 3,
or in other words the quantum matrix group for the quantum space 
$\cala$. 

The quantum matrix groups $\calh(\cala)$ for a quadratic 
$\cala$ are contained in the
Ewen and Ogievetsky classification \cite{EO} (see also\cite{O}) 
of the $GL(3)$ quantum matrix group. 
The quantum groups $\calh(\cala)$ for the cubic regular algebras
$\cala$ of dimension 3 are to the best of  our knowledge new ones 
(first reported by the author in \cite{TP}). 

One expects that all $\calh(\cala)$ have polynomial growth.
It is natural 
for the automorphism group $\calh(\cala)$ of 
an algebra $\cala$ with 
dimension 3 to expect $gk$-$\dim \calh(\cala)=9=3^2$. 
as it happens for the quadratic algebras.
Surprisingly one has  
$gk$-$\dim \calh(\cala)=7$
for some cubic $\cala$, e.g. of type  $S_{2}$.
 

%


\section*{Acknowledgments}
We would like to thank Dubois-Violette for many enlightning discussions.  
We thank   Roland Berger for his elucidating letter 
with critics of \cite{TP}. 
The Author was partially  supported by the  Bulgarian National Council
for Scientific Research under the project PH-1406, the  Euclid Network
HPRN-CT-2002-00325 and the European RTN
{\em ``Forces-Universe''} MRTN-CT-2004-005104.


\end{document}